\documentclass[a4paper,12pt]{amsart}

\usepackage{epsfig}
\usepackage{graphicx}
\usepackage[latin1]{inputenc}
\usepackage[T1]{fontenc}
\usepackage{indentfirst}
\usepackage{amssymb}
\usepackage{eufrak}
\usepackage{amsmath}
\usepackage{amsfonts}
\usepackage{amsthm}
\usepackage{mathrsfs}
\usepackage[bookmarks]{hyperref}
\usepackage{xypic}
\newcommand{\Log}{\mathop{\rm Log}\nolimits}
\newcommand{\val}{\mathop{\rm val}\nolimits}
\newcommand{\Val}{\mathop{\rm Val}\nolimits}
\newcommand{\ord}{\mathop{\rm ord}\nolimits}
\newcommand{\Arg}{\mathop{\rm Arg}\nolimits}

\newcommand{\Int}{\mathop{\rm Int}\nolimits}

\newcommand{\Sym}{\mathop{\rm Sym}\nolimits}

\newcommand{\supp}{\mathop{\rm supp}\nolimits}

\newcommand{\Critv}{\mathop{\rm Critv}\nolimits}
\newcommand{\Verte}{\mathop{\rm Vert}\nolimits}
\newcommand{\Mirr}{\mathop{\rm Mirr}\nolimits}

\input epsf

\def \square{\smallskip \hfill \vrule width 5 pt height 7
pt depth - 2 pt \smallskip }

\newenvironment{prooof}
{\noindent {{\it Proof} \;}}{\hspace*{\fill}\square\vskip 8pt}

\theoremstyle{plain}
\newtheorem{Lem}[subsection]{Lemma}
\newtheorem{The}[subsection]{Theorem}

\newtheorem{Pro}[subsection]{Proposition}
\theoremstyle{definition}
\newtheorem{Rem}[subsection]{Remark}

\newtheorem{Def}[subsection]{Definition}
\newtheorem{Exe}[subsection]{Example}
\newtheorem{Exes}[subsection]{Examples}

\begin{document}

\title[Complex tropical localization and mirror hypersurfaces]{Complex tropical localization, coamoebas, and  mirror tropical hypersurfaces}

\author{Mounir Nisse}
\address{Institut de Math{\'e}matiques de Jussieu (UMR 7586), Universit{\'e}
  Pierre et Marie Curie, Analyse Alg{\'e}brique\\ 175, rue du 
Chevaleret,\\ 75013
Paris} 
\email{nisse@math.jussieu.fr} 
\maketitle
%
%
%
%
%
\begin{abstract}
We introduce in this paper the concept of tropical mirror hypersurfaces and
we prove a complex tropical localization Theorem which is a  version of Kapranov's Theorem   \cite{K-00} in  tropical geometry. We give a geometric and a topological equivalence between  coamoebas of complex algebraic hypersurfaces defined by a maximally sparse polynomial and  coamoebas of maximally sparse complex tropical hypersurfaces. 
\end{abstract}

\setcounter{tocdepth}{1} \tableofcontents

\section{Introduction}

Amoebas have proved to be a very useful tool in several areas of mathematics, and  they have many applications in real algebraic geometry , complex analysis, mirror symmetry, algebraic statistics and in several other areas (see \cite{M1-02}, \cite{M2-04}, \cite{M3-00}, \cite{FPT-00}, \cite{PR1-04}, \cite{PS-04} and \cite{R-01}). They degenerate to a piecewise-linear object called {\em tropical varieties}, \cite{M1-02}, \cite{M2-04}, and \cite{PR1-04}.
 However, we can use amoebas as an intermediate link between the classical and the tropical geometry.  Coamoebas have a close relationship and similarities with amoebas and can be also used as an intermediate link between the tropical and the complex geometry. 

 A tropical hypersurface  is the set of points in $\mathbb{R}^n$ where some piecewise affine linear function (called tropical polynomial) is not differentiable. Such a
tropical polynomial may contain a tropical monomials which are not essential for the construction of the tropical hypersurface, but in the classical polynomial those monomials have a contribution and they often play a vital role in the geometry and the topology of the complex tropical hypersurface coamoeba defined by that polynomial.  In this paper, we give a process  to constructing the coamoeba of a complex tropical hypersurface by
 using a construction of a symmetric tropical hypersurface, which we call  a {\em mirror tropical hypersurface}, that allows us to see and to understand the contribution on the coamoeba of the non-essential monomials in the tropical polynomial. The construction consists to look at the deformation of the extending Newton polytope in $\mathbb{R}^n\times\mathbb{R}$ instead of the deformation of the tropical hypersurface itself. What is the same by duality, but in plus we have  a geometric point of view of that deformation.  A symmetry appears naturally in this deformation, whose center is the time when the dual subdivision $\tau$ of the Newton polytope $\Delta$ is reduced to one element (i.e., $\tau = \{\Delta \}$ itself). 

If $V_f$ is an algebraic hypersurface in $(\mathbb{K}^*)^n$ with $\mathbb{K}$ the field of the Puiseux series, then we obtain the following results:

\vspace{0.2cm}

\begin{The}[{\it Complex tropical localization}] Let $H_{\alpha\gamma}$ be  a hyperplane in $\mathbb{R}^n$ codual to an  edge $E_{\alpha\gamma}$ of the subdivision $\tau$, and $\mathcal{C}$ be a connected component of $H_{\alpha\gamma}\cap \Arg (V_{\infty ,\, f})$. Then we have one of the two following cases:
\begin{itemize}
\item[(i)] The dimension of $\mathcal{C}$ is $n-1$ and its interior is contained in the interior of a regular part of  $\Arg (V_{\infty ,\, f})$;
\item[(ii)] the dimension of  $\mathcal{C}$ is zero (i.e., discrete) and $\mathcal{C}$ is contained in the intersection of $H_{\alpha\gamma}$ and a line  codual to some proper face of $\Delta_v$.\end{itemize}
\end{The}

\begin{The} Let $V_f\subset (\mathbb{K}^*)^n$ be a hypersurface defined by a polynomial $f$ with Newton polytope $\Delta$ such  that the subdivision $\tau_f = \{ \Delta_1,\ldots ,\Delta_l \}$ dual to the tropical hypersurface $\Val (V_f)$ is a triangulation. Then the geometry and the topology of the complex tropical hypersurfaces $W(V_f)$  coamoebas are  completely determined and constructed  by  gluing  those of the truncated complex tropical hypersurfaces $W(V_{f^{\Delta_i}})$ using the complex tropical localization.
\end{The}

\vspace{0.2cm}

 If $V_f$ is a complex algebraic hypersurface, then we have the following result:

\vspace{0.2cm}

\begin{The} 
Let $V_f$ be a complex algebraic hypersurface defined by a maximally sparse polynomial $f$. Then there exist a  
 deformation of $V_f$ given by a  family of polynomials $f_t$ such that 
the coamoeba of the complex tropical hypersurface $V_{\infty ,\, f}$ (which is the limit of the $H_t(V_{f_t})$ with respect of the Hausdorff metric on compact sets of $(\mathbb{C}^*)^n$) has the same topology as
the coamoeba $co\mathscr{A}_{f}$ of $V_f$ (i.e., they are homeomorphic).
\end{The}

\vspace{0.2cm}

We recall the definitions and some Theorems of tropical geometry in section 2 alongside with all necessary notation. In section 3, we give the definition of complex tropical hypersurface and we describe those defined by maximally sparse polynomial with Newton polytope a simplex, and we give some examples of complex algebraic plane curves. In section 4, we introduce the notion of mirror tropical hypersurface, we give some examples, and we prove Theorem 1.2. In section 5, we prove the complex tropical localization Theorem. In section 6, we give a geometric and a topological description from the complex tropical hypersurface coamoeba to that of the complex algebraic hypersurface, and we will prove Theorem 1.3. Finally in section 7, we  give the geometric and topological description of the coamoebas of some complex algebraic plane curves.

\section{Preliminaries}

Let $\mathbb{K}$ be the field of the Puiseux series with real power, which is the field of the  series $\displaystyle{a(t) = \sum_{r\in A_a}\xi_rt^r}$ with $\xi_r\in \mathbb{C}^* =\mathbb{C}\setminus \{ 0\}$ and $A_a\subset \mathbb{R}$ is well-ordered set (which means that any subset has a smallest element); the smallest element of $A_a$ is called the {\em order} of $a$, and denoted by  $\ord (a) := \min A_a$. It is well known that the field $\mathbb{K}$ is algebraically closed and has a  characteristic equal to zero, and it has a non-Archimedean valuation $\val (a) = - \min A_a$ satisfying to the following properties.
\[ 
\left\{ \begin{array}{ccc}
\val (ab)&=& \val (a) + \val (b) \\
\val (a + b)& \leq& \max \{ \val (a)  ,\, \val (b)  \} ,
\end{array}
\right.
\]
and we put $\val (0) = -\infty$. If we denote by $\mathbb{K}^* =\mathbb{K}\setminus \{ 0\}$ and 
we apply the valuation map  coordinate-wise  we obtain a map $\Val : (\mathbb{K}^*)^n \rightarrow \mathbb{R} \cup \{ -\infty\}$ which we  will also call the valuation map.

\noindent If $ a \in \mathbb{K}^*$ is the Puiseux series
$\displaystyle{a = \sum_{j\in A_a}\xi_jt^j}$ with $\xi\in \mathbb{C}^*$
and $A_a\subset \mathbb{R}$ is a well-ordered set. We complexify the valuation
map as follows :
\[
\begin{array}{ccccl}
w&:&\mathbb{K}^*&\longrightarrow&\mathbb{C}^*\\
&&a&\longmapsto&w(a ) = e^{\val (a )+i\arg (\xi_{-\val
    (a )})}
\end{array}
\]
Let  $\Arg$ be the argument map $\mathbb{K}^*\rightarrow S^1$ defined by: for any $ a \in \mathbb{K}$ a Puiseux series so that
$\displaystyle{a = \sum_{j\in A_a}\xi_jt^j}$, then $\Arg (a) = e^{i\arg (\xi_{-\val (a)})}$ (this map extends the map $\Arg : \mathbb{C}^*\rightarrow S^1$ defined by $\rho e^{i\theta} \mapsto e^{i\theta}$).

Applying this map coordinate-wise we obtain a map :
\[
\begin{array}{ccccl}
W:&(\mathbb{K}^*)^n&\longrightarrow&(\mathbb{C}^*)^n
\end{array}
\]

\vspace{0.1cm}

\begin{Def}  The set  $V_{\infty}\subset (\mathbb{C}^*)^n$
  is a complex tropical hypersurface if and only if there
  exists an algebraic hypersurface
  $V_{\mathbb{K}}\subset(\mathbb{K}^*)^n$ over $\mathbb{K}$ such that
  $\overline{W(V_{\mathbb{K}})} = V_{\infty}$, where $\overline{W(V_{\mathbb{K}})}$ is the closure of $W(V_{\mathbb{K}})$ in  $(\mathbb{C}^*)^n \approx \mathbb{R}^n\times (S^1)^n$ as  a  Riemannian manifold with the metric of  the product of the
Euclidean metric on $\mathbb{R}^n$ and the flat metric on $(S^1)^n$.
\end{Def}

\vspace{0.2cm}

\noindent Let $V_f\subset (\mathbb{K}^*)^n$ be the algebraic hypersurface defined by the non-Archimedean polynomial:
$$
f(z) = \sum_{\alpha\in A} a_{\alpha} z^{\alpha}, \quad\quad\quad\quad\quad  z^{\alpha} = z_1^{\alpha_1}z_2^{\alpha_2}\ldots z_n^{\alpha_n}
$$
with $a_{\alpha}\in \mathbb{K}^*$ and $A$ a finite subset of $\mathbb{Z}^n$. We denote by $\Delta_f$ the Newton polytope of $f$, which is the convex hull in $\mathbb{R}^n$ of $A$. Let $\nu_f$ be the map defined on $A$ as follows:
\[
\begin{array}{ccccl}
\nu_f&:&A&\longrightarrow&\mathbb{R}\\
&&\alpha&\longmapsto&\ord (a_{\alpha}).
\end{array}
\]
\noindent The Legendre transform $\mathcal{L}(\nu_f)$ of the map $\nu_f$ is the piecewise affine linear convex function defined by:
\[
\begin{array}{ccccl}
\mathcal{L}(\nu_f)&:&\mathbb{R}^n&\longrightarrow&\mathbb{R}\\
&&x&\longmapsto&\max \{ <x,\alpha > - \nu_f(\alpha )\};
\end{array}
\]
where $<,>$ denotes the scalar product in the Euclidean space.

\begin{Def} The Legendre transform $\mathcal{L}(\nu_f)$ of the map $\nu_f$ is  called the {\em tropical polynomial} associated to $f$, and denoted by $f_{trop}$.
\end{Def}

\begin{The}[Kapranov, (2000)] The image of the algebraic hypersurface $V_f$ under the valuation map $\Val$ is the set $\Gamma_f$ of points in $\mathbb{R}^n$ where the piecewise affine linear function $f_{trop}$ is not differentiable.
\end{The}

\noindent We denote by $\tilde{\Delta}_f$ the extended Newton polytope of $f$ which is the convex hull of the subset $\{ (\alpha , \nu_f(\alpha )) \in A\times\mathbb{R}\}$ of $\mathbb{R}^n\times\mathbb{R}$. Let $\rho$ be the following map:
\[
\begin{array}{ccccl}
\rho&:&\Delta_f&\longrightarrow&\mathbb{R}\\
&&x&\longmapsto&\min \{ t\mid\, (x ,t)\in \tilde{\Delta}_f \}.
\end{array}
\]
It's clear that the linearity domains of $\rho$ define a convex subdivision $\tau_f = \{\Delta_1,\ldots ,\Delta_l\}$ of $\Delta_f$ (by taking the linear subsets of the lower boundary of $\tilde{\Delta}_f$, see \cite{PR1-04}, \cite{RST-05}, and \cite{IMS-07} for more details). Let $y= <x,v_i>+r_i$ be the equation of the hyperplane $H_i\subset \mathbb{R}^n\times\mathbb{R}$ containing the points with coordinates $(\alpha ,\nu_f(\alpha ))$ with $\alpha \in \Verte (\Delta_i)$.

\noindent There is a duality between the subdivision $\tau_f$ and the subdivision of $\mathbb{R}^n$ induced by $\Gamma_f$ (see \cite{PR1-04}, \cite{RST-05}, and \cite{IMS-07}), where each connected component of $\mathbb{R}^n\setminus \Gamma_f$ is dual to some vertex of $\tau_f$ and each $k$-cell of $\Gamma_f$ is dual to some $(n-k)$-cell of $\tau_f$. In particular, each $(n-1)$-cell of $\Gamma_f$ is dual to some edge of $\tau_f$.  If $x\in E_{\alpha\beta}^*\subset \Gamma_f$,  then $<\alpha , x> -\nu_f(\alpha ) = <\beta , x> -\nu_f(\beta )$, so $<\alpha  -\beta , x - v_i>  = 0$. This means that $v_i$ is a vertex of $\Gamma_f$ dual to some $\Delta_i$ having $E_{\alpha\beta}$ as edge.

\vspace{0.2cm}

Let $V$ be an algebraic hypersurface in $(\mathbb{C}^*)^n$ defined by
the complex polynomial:
$$
f(z) =\sum_{\alpha\in \supp (f)} a_{\alpha}z^{\alpha},\,\,\,\,\quad\quad\quad
z^{\alpha}=z_1^{\alpha_1}z_2^{\alpha_2}\ldots z_n^{\alpha_n},
$$
where $a_{\alpha}$ are non-zero complex numbers and $\supp (f)$ is
the support of $f$, and we
denote by $\Delta$ the Newton polytope of $f$ (i.e., the convex
hull in $\mathbb{R}^n$ of $\supp (f)$).

\vspace{0.1cm}

\noindent 
The following definition is given by M. Gelfand, M.M. Kapranov
 and A.V. Zelevinsky in \cite{GKZ-94}:
 \begin{Def} The {\em amoeba}  $\mathscr{A}$ of an algebraic hypersurface $V\subset
 (\mathbb{C}^*)^n$
 is  the image of $V$ under the map :
\[
\begin{array}{ccccl}
\Log&:&(\mathbb{C}^*)^n&\longrightarrow&\mathbb{R}^n\\
&&(z_1,\ldots ,z_n)&\longmapsto&(\log\parallel z_1\parallel ,\ldots ,\log\parallel
z_n\parallel ).
\end{array}
\]
\end{Def}

\vspace{0.1cm}

\noindent 
It was shown by M. Forsberg, M. Passare and A. Tsikh in \cite{FPT-00} that
there is an injective map between the set of components
$\{E_{\nu}\}$ of $\mathbb{R}^n\setminus \mathscr{A}$ and
$\mathbb{Z}^n\cap\Delta$:
$$
\ord :\{E_{\nu}\} \hookrightarrow \mathbb{Z}^n\cap\Delta
$$

\begin{The}[Foresberg-Passare-Tsikh, (2000)] Each  component
  of $\mathbb{R}^n\setminus \mathscr{A}$  is a convex domain and there
  exists a locally constant function:
$$
\ord :\mathbb{R}^n\setminus \mathscr{A} \longrightarrow \mathbb{Z}^n\cap\Delta
$$
which maps different components of the complement of $\mathscr{A}$
to different lattice points of $\Delta$.
\end{The}

\vspace{0.2cm}

 \noindent
The coordinates $z_j$ of
$z\in (\mathbb{C}^*)^n$ are parameterized by $z_j=\rho_j e^{i\arg
(z_j)}$ with $\rho_j =\,\parallel z_j\parallel \, \in ]0, \infty[$ and $\arg
(z_j)\in [0,2\pi[$ for $j=1,\ldots ,n$.
 Passare and Tsikh introduced the following set associated to a complex algebraic varieties.

\begin{Def}[Passare-Tsikh]  The {\em Coamoeba} $co\mathscr{A}\subset (S^1)^n$
of $f$ is the image of $V$ under the argument map $\Arg$ defined by the following:
\[
\begin{array}{ccccl}
\Arg&:&(\mathbb{C}^*)^n\approx\mathbb{R}^n\times (S^1)^n&\longrightarrow&(S^1)^n\\
&&(z_1,\ldots ,z_n)&\longmapsto&(e^{i\arg (z_1)},\ldots ,e^{i\arg (z_n)}).
\end{array}
\]
\end{Def}

\vspace{0.1cm} 

\section[]{Complex tropical hypersurfaces with a simplex Newton polytope}

\vspace{0.2cm}

Let $a=(a_1,\ldots ,a_n)\in (\mathbb{K}^*)^n $ and  $H_a\subset(\mathbb{K}^*)^n$ be the hyperplane defined by the polynomial $f_a(z_1,\ldots ,z_n) = 1+\sum_{j=1}^na_jz_j$, then it's clear that $H_a = \tau_{a^{-1}}(H_{1})$. Let $L$ be an invertible matrix with integer coefficients and positive determinant
$$
L = \left(
\begin{array}{ccc}
\alpha_{11}&\hdots&\alpha_{1n}\\
\vdots&\ddots&\vdots\\
\alpha_{n1}&\hdots&\alpha_{nn}
\end{array}
\right) ,
$$
and let $\Phi_{L,\, a}$ be the homomorphism of the algebraic torus defined as follow.
\[
\begin{array}{ccccl}
\Phi_{L,\, a}&:&(\mathbb{K}^*)^n&\longrightarrow&(\mathbb{K}^*)^n\\
&&(z_1,\ldots ,z_n)&\longmapsto&(a_1\prod_{j=1}^nz_j^{\alpha_{j1}},\ldots ,a_n\prod_{j=1}^nz_j^{\alpha_{jn}}).
\end{array}
\]
Let $V_f\subset (\mathbb{K}^*)^n$ be the hypersurface defined by the polynomial 
$$
f(z_1,\ldots ,z_n)=1+\sum_{k=1}^na_k\prod_{j=1}^nz_j^{\alpha_{jk}},
$$
such that its Newton polytope is the simplex $\Delta_f$ that is the  image by $L$ of the standard simplex. The matrix $L$ is invertible, so $\Phi_{L,\, a}(V_f)=H_a$, and then  ${}^tL^{-1}(\Val (H_a)) = \Val (V_f)$. It was the same thing  for the complex tropical hypersurface i.e., ${}^tL^{-1}(W(V_f))=W(H_a)$,  (because for any $k=1,\ldots ,n$ we have  $\arg (a_k\prod_{j=1}^nz_j^{\alpha_{jk}})= \arg (a_k)+ \sum_{j=1}^n<\alpha_{jk}, \arg (z_j)>$), abuse of notations; to be more precise we have ${}^tL^{-1}(\Log (\Arg (W(V_f)))) = \Log (\Arg (W(H_a)))$. Hence we have the following (for more details, see \cite{N1-07}): 
$$
co\mathscr{A}(V_f) = \tau_{{}^tL^{-1}(a^{-1})}\circ {}^tL^{-1} (co\mathscr{A}(H_1)).
$$
So, the coamoeba of any hypersurface defined by a maximally sparse polynomial (that the number of its coefficients is equal to the number of its Newton polytope vertices) with a simplex as Newton polytope,  can be easily drawn.
We  remark that the field of Puiseux series $\mathbb{K}$ can be replaced by the field of complex numbers and we have the same results with the same formulas.

\vspace{0.1cm}

\begin{Exe}

We draw on figure 1 the coamoeba  of the complex  curve defined by the polynomial $f_1(z,w) = w^3z^2+wz^3+1$
where the matrix  ${}^tL_1^{-1}$ is equal to $\frac{1}{7}\left(
  \begin{array}{cc} 3&-1\\ -2&3\end{array}\right) $
and on figure 2 the coamoeba of the complex curve 
defined by the polynomial $f_2(z,w) = w^2z^2+z+w$
where the matrix  ${}^tL_2^{-1}$ is equal to $\frac{1}{3}\left(
  \begin{array}{cc} 1&1\\ -2&1\end{array}\right)$.

\begin{figure}[h!]
\begin{center}
\includegraphics[angle=0,width=0.4\textwidth]{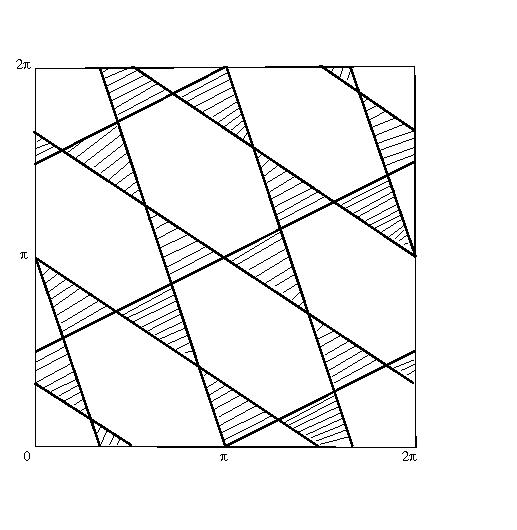}
\caption{The coamoeba of the curve defined by the polynomial $f_1(z,w)=wz^3+z^2w^3+1$}
\label{c}
\end{center}
\end{figure}

\begin{figure}[h!]
\begin{center}
\includegraphics[angle=0,width=0.4\textwidth]{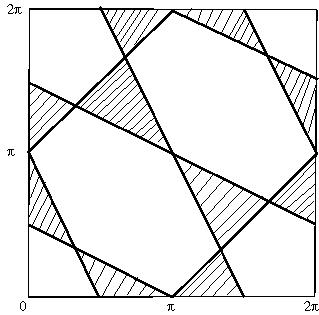}
\caption{The coamoeba of the curve defined by the polynomial $f(z,w)=z+w+z^2w^2$}
\label{c}
\end{center}
\end{figure}

\end{Exe}

\vspace{0.2cm}

\section{Tropical mirror hypersurfaces}

\vspace{0.2cm}

Let $V_f\subset (\mathbb{K}^*)^n$ be an algebraic hypersurface defined by a polynomial $f$ and   we 
assume  that $\Delta_f$ is a simplex, $A = \Verte (\Delta_f)\cup \{\beta\}$, and the coefficient $a_{\beta}$ is monomial. In addition, we suppose that $\beta\in \Verte (\tau_f)$. Let $\displaystyle{\{f_u\}_{u\in]-1;1]}}$ be the family of polynomials defined as follow:
$$
f_u(z) = a_{\beta ,\, u}z^{\beta} + \sum_{\alpha\in\Verte (\Delta_f)} a_{\alpha}z^{\alpha},
$$
with $a_{\beta ,\, u}$ such that:
\[
a_{\beta ,\, u} = \left\{
\begin{array}{ll}
\xi_{\beta} t^{u\nu_f(\beta ) + (1-u)(<\beta ,v> + r)}& \mbox{if\, $u\in [0;1]$}\\
\xi_{\beta} t^{(1-u)(<\beta ,v> + r) -\frac{u}{u+1}}&\mbox{if\, $u\in ]-1;0]$}
\end{array}
\right.\]
where $\xi_{\beta}$ is the complex coefficient of $a_{\beta}$, and $v$ is the vertex of the tropical hypersurface $\Val (V_f)$. We can assume that $<\beta ,v> + r \geq 0$ (multiplying $f$ by a power of $t$ if necessary); so the map  $u\longmapsto (1-u)(<\beta ,v> + r) - \frac{u}{u+1}$ is a decreasing function  on $]-1; 0]$.


\begin{Rem}${}$

\begin{itemize}
\item[(a)]\, The deformation given above is such that  $f_1(z) = f(z)$;
\item[(b)]\, for any $u\in ]-1;0]$, the subdivision $\tau_{f_u} =\{ \Delta_f\}$ (i.e., trivial);
\item[(c)]\, with the same assumption as above, when the order of the monomial $a_{\beta ,\, u}$ reaches the hyperplane of $\mathbb{R}^n\times\mathbb{R}$ containing the points with coordinates $(\alpha , \nu_f(\alpha ))$ and $\alpha\in \Verte (\Delta_f)$ (i.e., for $u\leq 0$). Then we consider the family of polynomials $\tilde{f}_u$ defined by:
$$ 
\tilde{f}_u(z) = \tilde{a}_{\beta ,\, u} z^{-\beta} + \sum_{\alpha\in\Verte (\Delta_f)}\tilde{a}_{\alpha}z^{-\alpha},
$$
such that if \, $\alpha\in \Verte (\Delta_f)$ and 
 $\displaystyle{a_{\alpha}(t)=\sum_{r\geq \ord (a_{\alpha})}\xi_{\alpha , r}t^r}$, then  we set  $\displaystyle{\tilde{a}_{\alpha}(t)=}$\\
$\displaystyle{\sum_{r\geq \ord (a_{\alpha})}\xi_{\alpha , r}t^{-r}}$
 and if\,  $a_{\beta ,\, u}(t)=\xi_{\beta} t^{\ord (a_{\beta ,\, u})}$, we set \, 
$\tilde{a}_{\beta ,\, u}(t)= \xi_{\beta} t^{-\ord (a_{\beta ,\, u})}$. In this case, we have  the convergence  when $t$ tends to the infinity, because the induced  transformation of $\mathbb{K}$ is given by $t\longmapsto t^{-1}$.
Let  $\mathscr{I} : (\mathbb{K}^*)^n\rightarrow (\mathbb{K}^*)^n$ be the transformation defined as  $(z_1,\ldots ,z_n)\mapsto (z_1^{-1},\ldots ,z_n^{-1})$, then by 
making the change of the variable $t = \frac{1}{\tau}$, we can see that $\tilde{f}_u(z)= f_u\circ\mathscr{I}(z)$, and then $V_{\tilde{f}_u} = V_{f_u\circ\mathscr{I}} = \mathscr{I}(V_{f_u})$. The tropical polynomial associated to $\tilde{f}_u$ is given by 
$$
\tilde{f}_{u,\, trop} = \max_{\alpha\in \Sym (A)}\{ <x, \gamma> -\val (a_{\gamma , u})\}
$$
with $\Sym (A)$ the subset of $\mathbb{Z}^n$ symmetric to $A$ relatively to the origin. There exist a positive number $s\in ]-1;0]$
such that the non-Archimedean amoebas defined by the tropical polynomials $\tilde{f}_{u,\, trop}$ with $u\in [-s;0]$ are symmetric to those defined by $f_{u,\, trop}$ with $u\in [0;1]$ (By an automorphism of 
 $(\mathbb{K}^*)^n$ if necessary, we can assume that $\val (a_{\alpha})=0$ for any $\alpha\in \Verte (\Delta_f$, and in this case $s=-\nu_f(a_{\beta })$.). 
So, 
we can apply now  Kapranov's theorem to the tropical hypersurfaces $\Gamma_{\tilde{f}_{u}}$, and from  the equality $V_{\tilde{f}_u} = \mathscr{I}(V_{f_u})$, we deduced that the coamoeba of $V_{\infty ,\, f_u}=W(V_{f_u})$ is the symmetric of the coamoeba of $V_{\infty ,\, \tilde{f}_u}=W(V_{\tilde{f}_u})$.
\item[(d)]\,  One way to look at the deformation of a tropical hypersurface is to think of it as a deformation of the extending Newton polytope of its defining polynomial. More precisely,
the deformation $\tilde{f}_{u,\, trop}$ can be seen as a continuation of the deformation of the normal vectors to the hyperplanes in $\mathbb{R}^n\times\mathbb{R}$ containing the lifting of the $\Delta_i$'s element of the subdivision $\tau_{f_u}$ dual to ${f}_{u,\, trop}$ with $0\leq u \leq 1$. Indeed, when $u=0$, all the normal vectors are equal and then for $u\leq 0$ the coefficient of index $\beta$ becomes inessential in the tropical polynomial  ${f}_{u,\, trop}$, but in the non-Archimedean polynomial $f_u$, it has a contribution and plays a crucial role for the determination  of the complex tropical hypersurface coamoeba.

\end{itemize}
\end{Rem}

\begin{Def} The tropical hypersurfaces $\Val (V_{\tilde{f}_{u}})$ 
defined by the tropical polynomials $\tilde{f}_{u,\, trop}$ for $u\in ]-1;0]$, are called the  {\em tropical mirror} for  the hypersurfaces $V_{f_{u}}$  defined by the  polynomial $f_{u}$. We denote this hypersurface by $\Mirr_{trop} (V_{f_{u}}):=\Val (V_{\tilde{f}_{u}})$.
\end{Def}

\noindent We can see that if $\Gamma$ is a  tropical hypersurface  with only one vertex, and $V_1, V_2$ are two hypersurfaces  in $(\mathbb{K}^*)^n$ such that $\Val (V_i)=\Gamma$ for $i=1,2$, then the two mirror tropical hypersurfaces $\Mirr_{trop} (V_1)$ and $\Mirr_{trop} (V_2)$ are not necessary the same. A similar algebraic construction is given by Z. Izhakian and L. Rowen in \cite{IR-08}.

\begin{The} Let $V_f\subset (\mathbb{K}^*)^n$ be a hypersurface defined by a polynomial $f$ with Newton polytope $\Delta$, and assume that the subdivision $\tau_f = \{ \Delta_1,\ldots ,\Delta_l \}$ dual to the tropical hypersurface $\Val (V_f)$ is a triangulation. Then the geometry and the topology of the complex tropical hypersurfaces $W(V_f)$  coamoebas are completely determined and constructed  by  gluing  those of the truncated complex tropical hypersurfaces $W(V_{f^{\Delta_i}})$ using the complex tropical localization.
\end{The}

\begin{prooof}
Suppose that $V_f\subset (\mathbb{K}^*)^n$ is defined by the polynomial $\displaystyle{f(z) =\sum_{\alpha\in A} a_{\alpha}z^{\alpha}}$ with Newton polytope $\Delta$ equal to the convex hull of $A$, and let $A_i=A\cap\Delta_i$. If $f^{\Delta_i}$ denotes the truncation of $f$ to $\Delta_i$, then the assumption of Theorem 4.3, means that the spine of the hypersurface amoeba of $V_{f^{\Delta_i}}$ has only one vertex. Let $\tau_i = \cup_j \Delta_{ij}$ be the convex subdivision of $\Delta_i$ given by taking the upper bound of the convex hull of the set $\{ (\alpha ,r)\in A_i\times \mathbb{R} \mid r\leq \ord (a_{\alpha})\}$, which we can suppose to be  a triangulation (by a small perturbation of the coefficients order if necessary). Let $inv : \mathbb{C}((t))\rightarrow \mathbb{C}((\rho ))$ be the morphism sending $t$ of valuation $+1$ to $\rho$ of valuation $-1$, and 
let $\tilde{f}$ be the polynomial defined by $\displaystyle{\tilde{f}(z)= \sum_{\alpha\in A_i} \tilde{a}_{\alpha}z^{-\alpha}}$ with $\tilde{a}_{\alpha}=inv(a_{\alpha})$ (this means that if $\displaystyle{a_{\alpha}(t)=\sum_{r\geq \ord (a_{\alpha})}\xi_{\alpha , r}t^r}$ then $\displaystyle{\tilde{a}_{\alpha}(t)= \sum_{r\geq \ord (a_{\alpha})}\xi_{\alpha , r}t^{-r}}$). We use now induction on the volume of $\Delta$, and we assume that the coamoeba of any $V_{f^{\Delta_{ij}}}$ is constructed for each index $ij$ using the complex tropical localization which we develop in the next section. By construction we have $V_{\tilde{f}_i}=V_{f^{\Delta_i}\circ\mathscr{I}}=\mathscr{I}(V_{f^{\Delta_i}})$. The coamoeba of $V_{\tilde{f}_i}$ can be constructed, because in this case, one can apply  Kapranov's Theorem, and we can also build the coamoeba of $V_{f^{\Delta_i}}$. Knowing now all the coamoebas of the $V_{f^{\Delta_i}}$'s, the coamoeba of the hypersurface $V_f$ itself
can be built by reusing Kapranov's Theorem.
\end{prooof}

\begin{Exes}${}$
\begin{itemize}
\item[(a)] Example of the parabola (see figures 3, 4, and 5), where the deformation is seen  as a  deformation of the normal vectors to the hyperplanes in $\mathbb{R}^n\times\mathbb{R}$ containing the lifting of the $\Delta_i$'s, and the points with coordinates $(\alpha ,\ord (\alpha ))$ and $\alpha\in\Verte (\Delta_f )$ are fixed.
\begin{figure}[h!]
\begin{center}
\includegraphics[angle=0,width=0.2\textwidth]{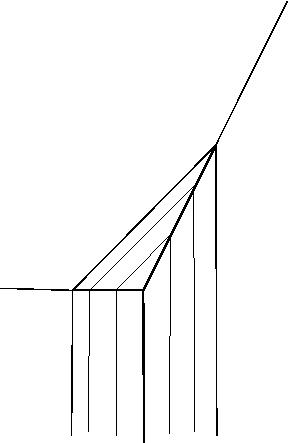}\quad\quad\quad
\includegraphics[angle=0,width=0.3\textwidth]{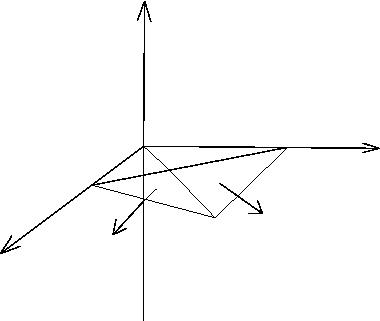}
\caption{The tropical curves $\Val (V_{f_u})$ for $u\in [0;1]$.}
\label{c}
\end{center}
\end{figure}
\begin{figure}[h!]
\begin{center}
\includegraphics[angle=0,width=0.2\textwidth]{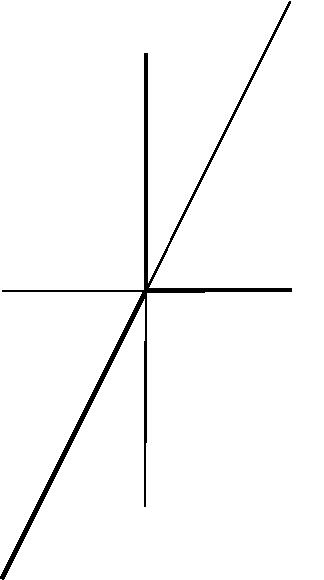}\quad\quad\quad
\includegraphics[angle=0,width=0.3\textwidth]{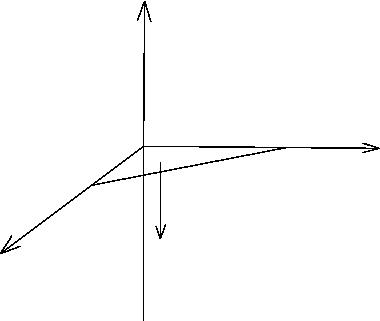}
\caption{The tropical curves $\Val (V_{f_0})$ and  $\Val (V_{\tilde{f}_0})$.}
\label{c}
\end{center}
\end{figure}
\begin{figure}[h!]
\begin{center}
\includegraphics[angle=0,width=0.2\textwidth]{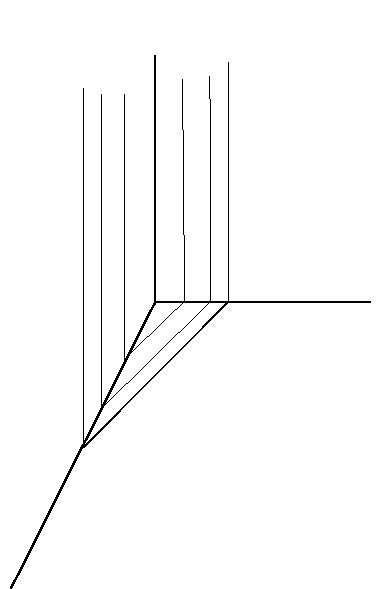}\quad\quad\quad
\includegraphics[angle=0,width=0.3\textwidth]{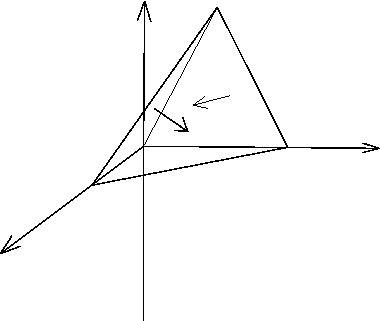}
\caption{The tropical curves $\Val (V_{\tilde{f}_u})$ for $u\in ]-1;0]$.}
\label{c}
\end{center}
\end{figure}

\item[(b)] We give here an example where $\beta\in \Int (\Delta_f)$ (see figures 6, 7, and 8), and as in the previous  example, the deformation is supposed to fix the order of the coefficients of index in the vertices of the Newton polygon.
\begin{figure}[h!]
\begin{center}
\includegraphics[angle=0,width=0.3\textwidth]{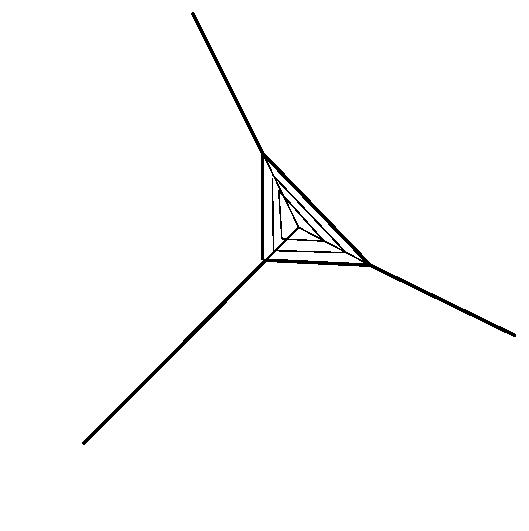}\quad\quad\quad\quad
\includegraphics[angle=0,width=0.2\textwidth]{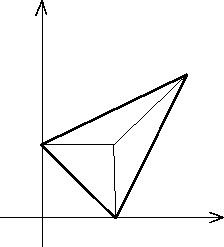}
\caption{The tropical curves $\Val (V_{f_u})$ for $u\in [0;1]$.}
\label{c}
\end{center}
\end{figure}
\begin{figure}[h!]
\begin{center}
\includegraphics[angle=0,width=0.3\textwidth]{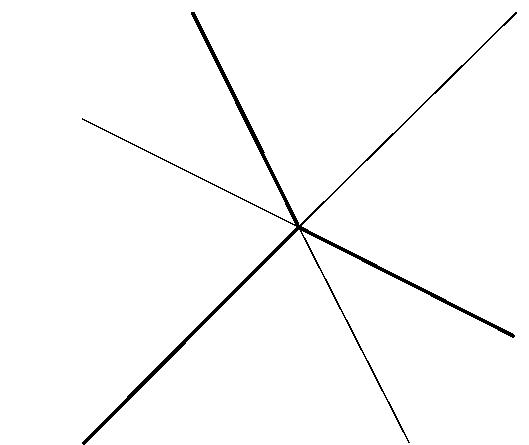}\quad\quad\quad\quad
\includegraphics[angle=0,width=0.3\textwidth]{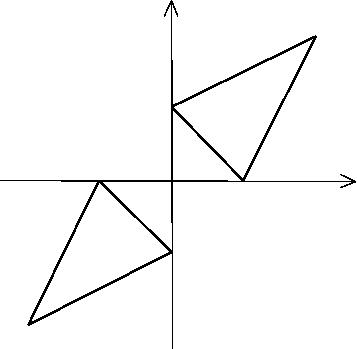}
\caption{The tropical curves $\Val (V_{f_0})$ and  $\Val (V_{\tilde{f}_0})$.}
\label{c}
\end{center}
\end{figure}
\begin{figure}[h!]
\begin{center}
\includegraphics[angle=0,width=0.2\textwidth]{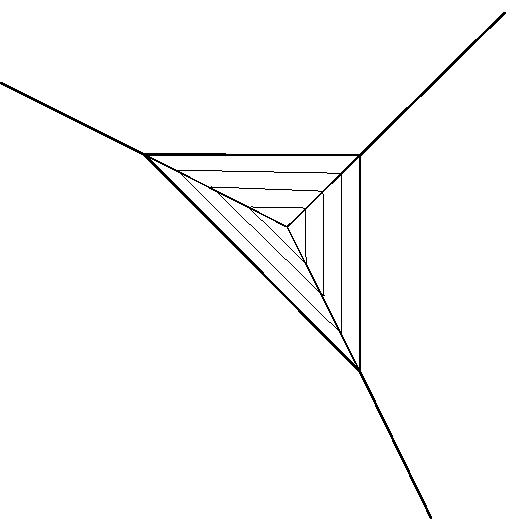}\quad\quad\quad\quad
\includegraphics[angle=0,width=0.2\textwidth]{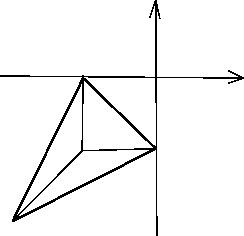}
\caption{The tropical curves $\Val (V_{\tilde{f}_u})$ for $u\in ]-1;0]$.}
\label{c}
\end{center}
\end{figure}
\end{itemize}
\end{Exes}

%
%
%
%
%
%
%
%



\section{Camoebas of complex tropical hypersurfaces}

\vspace{0.2cm}

In this section we consider an algebraic hypersurface $V$ over the field of Puiseux series $\mathbb{K}$ defined by a polynomial $f$ with Newton polytope $\Delta$. We denote by $\Gamma$ the non-Archimedean amoeba of $V$ and by $V_{\infty ,\, f}$ the complex tropical hypersurface image of $V$ under the map $W$. Let us denote by $\tau$ the subdivision of $\Delta$ dual to $\Gamma$ which we suppose to be  a triangulation, and assume that $f$ is defined  as follows:
$$
f(z) =\sum_{\alpha\in \supp (f)} a_{\alpha}z^{\alpha},\,\,\,\,\quad\quad\quad
z^{\alpha}=z_1^{\alpha_1}z_2^{\alpha_2}\ldots z_n^{\alpha_n},
$$
where $a_{\alpha}$ are non-zero complex Puiseux series and $\supp (f)$ is
the support of $f$.


\newpage

\subsection{Complex tropical localization} ${}$

\vspace{0.2cm}

\begin{Def}Let $\mathbb{R}^n$ be the universal covering of the real torus $(S^1)^n$. Let $\alpha$ and $\beta$ be  in the support of $f$. A hypersurface $H_{\alpha\beta}\subset \mathbb{R}^n$ is called {\em codual} (or corresponding) to an edge $E_{\alpha\beta}$ in $\tau$ if it is given  by the following equation:
$$
\arg (a_{\alpha}) - \arg (a_{\beta}) + <\alpha - \beta , x> = \pi .
$$
In addition if $E_{\alpha\beta}$ is an external edge of $\tau$ (i.e., $E_{\alpha\beta}$ is a proper edge of the Newton polytope $\Delta$), then  $H_{\alpha\beta}$ is called an {\em external hyperplane}.
\end{Def}

\begin{Def} An open  subset $\mathcal{C}$ of the coamoeba of a complex tropical hypersurface  $V_{\infty ,\, f}$  is called  {\em regular} if for any point $x$ in $\mathcal{C}$ there exist an open subset $V(x)$ in  $(S^1)^n$ containing $x$ with $V(x)\subset\mathcal{C}$  and an open subset $U$  in  $\mathbb{R}^n$ such that $V(x)\subset \Arg (\Log^{-1}(U)\cap V_{\infty ,\, f})_i$ where $(\Log^{-1}(U)\cap V_{\infty ,\, f})_i$ is one connected component of $\Log^{-1}(U)\cap V_{\infty ,\, f}$.

\end{Def}

\noindent We denote by $\Critv (\Arg )$ the set of critical values points in the coamoeba $co\mathscr{A}$ of a complex tropical hypersurface $V$.

\begin{Def} An {\em extra-piece} is a connected component $\mathcal{C}$ of $co\mathscr{A}\setminus \Critv (\Arg )$ such that the boundary of its closure $\partial \, \overline{\mathcal{C}}$ is not contained in the union of hyperplanes codual to the edges of the subdivision.
\end{Def}

 \noindent This means that its boundary contains at least one component (smooth) in the set of critical values of the argument map. In the following Lemma we assume that the subdivision $\tau$ of the Newton polytope $\Delta$ dual to the non-Archimedean amoeba $\Gamma$ is a triangulation and contains inner edges.

\vspace{0.2cm}

\noindent We begin by proving the following Lemma which is a local version of the  Theorem 6.1 in the complex tropical case.

\begin{Lem} Let $H_{\alpha\gamma}$ be  a hyperplane in $\mathbb{R}^n$ codual to an inner edge $E_{\alpha\gamma}$ of the subdivision $\tau$. Then any connected component $\mathcal{C}$ of $H_{\alpha\gamma}\cap \Arg (V_{\infty ,\, f})$ has a dimension $n-1$ and its interior is contained in the interior of a regular part of $\Arg (V_{\infty ,\, f})$.
\end{Lem}

\begin{prooof} Let $\Delta_{1}$ and $\Delta_{2}$ be  two elements of $\tau$ with a common edge $E_{\alpha\gamma}$, and  $v_{1}$ and $v_{2}$ be their dual vertices in  the non-Archimedean amoeba $\Gamma$. Let $\{ x_m\}$ be a sequence in $\Arg (\Log^{-1}(v_{1})\cap V_{\infty ,\, f^{\Delta_{1}}})\setminus \Arg (\Log^{-1}(v_{2})\cap V_{\infty ,\, f^{\Delta_{2}}})$ which converge to some point $x$ in $H_{\alpha\gamma}\setminus \Arg (\Log^{-1}(v_{2})\cap V_{\infty ,\, f^{\Delta_{2}}})$. Let $\mathscr{C}$ be a connected component of $\Arg^{-1}(\{ x_m\} ) \cap V_{\infty ,\, f}$ and $\{ z_m\}\subset \mathscr{C}$ be a sequence such that $\Arg (z_m) = x_m$ for each $m$. We claim that the sequence $\{ z_m\}$ (by taking a subsequence if necessary)  converges to some point $z$ in $V_{\infty ,\, f}$. Indeed, the sequence $\{ \Log (z_m)\}$ converge to $v_{2}$ because the argument of $z_m$ is $x_m$ which converges to $x\in H_{\alpha\gamma}$, and $x$ is  an infinite point for $\Arg (V_{\infty ,\, f^{\Delta_{1}}})$. This means that  $\{ \Log (z_m)\}$ converges asymptotically in  the direction of $E_{\alpha\gamma}$ to the infinity of $\Log (V_{\infty ,\, f^{\Delta_{1}}})$.  So $z_m$ converge to the point $z$ of $(\mathbb{C}^*)^n$ with argument $x$ and the valuation $v_{2}$. $V_{\infty ,\, f}$ is closed, hence  $z\in V_{\infty ,\, f}$. Then all the components of $H_{\alpha\gamma}\setminus \left( \overline{\Arg (\Log^{-1}(v_{1})\cap V_{\infty ,\, f^{\Delta_{1}}})}\cap \overline{\Arg (\Log^{-1}(v_{2})\cap V_{\infty ,\, f^{\Delta_{2}}})}\right)$
 are in the interior of $\Arg (V_{\infty ,\, f})$.  Let now $x$ be a point in the interior of the following set:
$$
H_{\alpha\gamma}\cap 
\left(
\overline{\Arg (\Log^{-1}(v_{1})\cap V_{\infty ,\, f^{\Delta_{1}}})}
\cap
\overline{\Arg (\Log^{-1}(v_{2})\cap \Arg (V_{\infty ,\, f^{\Delta_{2}}})}
\right) ,
$$
and $\{ x_m\}$ be a sequence in $\Arg (\Log^{-1}(v_{1})\cap V_{\infty ,\, f^{\Delta_{1}}})\cap \Arg (\Log^{-1}(v_{2})\cap \Arg (V_{\infty ,\, f^{\Delta_{2}}})$ such that $x_m$ converges to $x$. We claim that  there is no sequence $\{ z_m\}$ in $V_{\infty ,\, f}$ such that $\Arg (z_m) = x_m$ for any $m$ and $z_m$ converges in $V_{\infty ,\, f}$ to some point $z$ such that $\Arg (z) = x$. Indeed, assume on the contrary that there exists a sequence $\{ z_m\}$ in $V_{\infty ,\, f}$ satisfying the  assumption and converging to $z$ in $V_{\infty ,\, f}$. On one  hand we know that $\Log (z_m)$ converges to $v_{2}$, because the argument of $z_m$ converges to $x\in H_{\alpha\gamma}$ which is an infinite point for $\Arg (\Log^{-1}(v_{1})\cap V_{\infty ,\, f^{\Delta_{1}}})$ and then the valuation of the $z_m$'s tends to the infinity asymptotically in the direction of $E_{\alpha\gamma}$ to $v_{2}$ (because $v_{2}$ represents the infinity for $\Log (V_{\infty ,\, f^{\Delta_{1}}})$ in the direction of $E_{\alpha\gamma}$). On the other hand, for the same reasons, the sequence $\Log (z_m)$ converge to $v_{1}$. Contradiction, because by assumption $v_{1} \ne v_{2}$.
In this case we have the so-called extra-piece.
\end{prooof}

\begin{Pro} Let $H_{\alpha\gamma}$ be  a hyperplane in $\mathbb{R}^n$ codual to an external edge $E_{\alpha\gamma}$ of the subdivision $\tau$, and let  $\mathcal{C}$ be a connected component of $H_{\alpha\gamma}\cap \Arg (V_{\infty ,\, f})$. Then we have one of the two following cases:
\begin{itemize}
\item[(i)] The dimension of $\mathcal{C}$ is $n-1$ and its interior is contained in the interior of a regular part of  $\Arg (V_{\infty ,\, f})$;
\item[(ii)] the dimension of  $\mathcal{C}$ is zero (i.e., discrete) and $\mathcal{C}$ is contained in the intersection of $H_{\alpha\gamma}$ and a line  codual to some proper face of $\Delta_v$.\end{itemize}
\end{Pro}

 If the edge $E_{\alpha\gamma}$ is a common edge to more than one element of the subdivision $\tau$ (which can occur only if $n>2$), then by Lemma 5.5  we have  the first case. Assume that $E_{\alpha\gamma}$ is an edge of only one element $\Delta_v$ of $\tau$, and we denote by $v$ the vertex of the tropical hypersurface dual to $\Delta_v$. Let $z\in \Log^{-1}(v)\cap V_{\infty ,\, f}$ such that $\Arg (z) = x$ which we assume in $H_{\alpha\gamma}$. We denote by $\mathcal{C}$ the connected component of $H_{\alpha\gamma}\cap \Arg (V_{\infty ,\, f})$ containing $x$. We have to consider the following cases:
\begin{itemize}
\item[(a)]\, $\supp (f) = \Verte (\Delta_v)$, in this case there is nothing to prove, and we have case (ii) of the Proposition.
\item[(b)]\, $\supp (f) \cap \Delta_v = \Verte (\Delta_v)$ or  $\supp (f)  =  \Verte (\Delta_v) \cup \{ \beta_1,\ldots ,\beta_l \}$ with $\beta_j\in \Delta_v\cap \mathbb{Z}^n$ for any $j$.
\end{itemize}
All other cases will be easily deduced thereof.
Assume that $\supp (f)  =  \Verte (\Delta_v) \cup \{ \beta\}$ with $\beta\in \Delta_v$.

\vspace{0.1cm}

\begin{Lem} With the above notations, let $A$ be the interior of $\mathcal{C}$, then for any $x\in A$ there exists an open neighborhood $\mathscr{V}(x)$ of $x$ in $(S^1)^n$ such that $\mathscr{V}(x)\subset \Arg (V_{\infty ,\, f})$.
\end{Lem}

\begin{prooof}
\vspace{0.1cm}
Indeed, assume on the contrary that there exists a small open neighborhood $\mathscr{V}(x)$ of $x$ in $\mathbb{R}^n$ such that $\mathscr{V}(x)\cap co\mathscr{A}_{f^{\Delta_{\alpha_j}}}$ is empty, where $\Delta_{\alpha_j}$ is the simplex with vertices $\{ \alpha_1,\ldots ,\widehat{\alpha_j},\ldots ,\alpha_{n+1},\beta\}$ and $E_{\alpha\gamma} = E_{\alpha_1\alpha_2}$ with $j\ne 1, 2$ (here we use the same letter for
 $x$ and its lifting to the universal covering of the torus; abuse of notation). This means that $\mathscr{V}(x)\cap \Arg (V_{\infty ,\, f})$ lies in one side of the hyperplane $H_{\alpha_1\alpha_2}$. So the dominating monomials in $W^{-1}(\Arg^{-1}(\mathscr{V}(x)) \cap V_{\infty ,\, f})$ are $a_{\alpha_1},\widehat{a_{\alpha_2}},\ldots ,a_{\alpha_{n+1}, \beta}$, because if the monomial $a_{\alpha_2}$ is a dominating one, then $\Arg (V_{\infty ,\, f}) \cap \mathscr{V}(x)$ lies on both sides of $H_{\alpha_1\alpha_2}$. From  Remarks 4.1 (a), (b), (c) and  Kapranov's Theorem \cite{K-00}, we obtain that the dominating monomials in $W^{-1}(\Log^{-1}(v)\cap V_{\infty ,\, f})$ are $a_{\alpha_1},a_{\alpha_2},\ldots ,a_{\alpha_{n+1}}$. Hence $z$ lies in the  domain where the monomials $a_{\alpha_1},\widehat{a_{\alpha_2}},
\ldots ,a_{\alpha_{n+1}}$ are dominating (a proper face of the simplex $\Delta_v$),  and then $\Arg (z) = x$ is contained in $H_{\alpha_1\alpha_2}\cap \Arg (V_{\infty ,\, f^{\Delta_{\alpha_2}}})$. Contradiction, because $H_{\alpha_1\alpha_2}\cap \Arg (V_{\infty ,\, f^{\Delta_{\alpha_2}}})$ is discrete  and then the intersection of any open neighborhood of $x$ in $\mathbb{R}^n$ with $\Arg (V_{\infty ,\, f})$ lies on both  sides of the hyperplane $H_{\alpha_1\alpha_2}$. In this case we have some extra-piece.
 
\end{prooof}

Theorem 1.1 is an immediate consequence of Lemma 5.5 and  Proposition 5.6.

\vspace{0.2cm}

\section{Coamoebas of complex algebraic  hypersurfaces}

\vspace{0.2cm}

We now turn our attention to complex algebraic hypersurfaces, so in this section 
we assume that the polynomial $f$ is complex. We will give a caracterization of the argument map critical values set contained in the hyperplanes codual to the edges of the subdivision $\tau$ dual to the spine of the amoeba $\mathscr{A}$ of  $V_f$, and we have the following.

\begin{The} Let $H_{\alpha\gamma}$ be a hyperplane in $\mathbb{R}^n$ codual to an edge $E_{\alpha\gamma}$. Then the intersection $H_{\alpha\gamma}\cap \Critv (\Arg )$ is discrete and it is  contained in the union of lines $L_{\alpha '\beta '}$ codual to some faces of $\tau$.
\end{The}

\begin{prooof}
Assume that there is an open subset $A$ of $H_{\alpha\gamma}$ such that $A\subset co\mathscr{A}_{V_f}$, then we claim that $A\subset co\mathscr{A}_{V_{\infty ,\, f}}$.
%
Indeed, assume that $E_{\alpha\gamma}$ is a common edge for two simplices $\Delta_1$ and $\Delta_2$. Let $y= <x,a_1> + b_1$ be the equation of the hyperplane in $\mathbb{R}^n\times\mathbb{R}$ containing the points with coordinates $(\alpha ,\nu (\alpha ))$ and  $\alpha\in\Verte (\Delta_1)$ and $\nu$ the Passare-Rullg\aa rd function. Let $\displaystyle{f_t(z) = \sum a_{\alpha} (et)^{<\alpha , a_1> + b_1} z^{\alpha}} = (et)^{b_1}\sum a_{\alpha} ((et)^{a_1}z)^{\alpha}$ with $a_1 = (a_{11},\ldots ,a_{1n})$  and  
$$
((et)^{a_1}z)^{\alpha} = (et)^{a_{11}{\alpha_1}}z_1^{\alpha_1}(et)^{a_{12}{\alpha_2}}z_2^{\alpha_2}\ldots (et)^{a_{1n}{\alpha_n}}z_n^{\alpha_n}.
$$ 
Hence $V_{f_t}\subset (\mathbb{C}^*)^n$
 is the image of $V_f$ under the self diffeomorphism $\phi_{t}$ of $(\mathbb{C}^*)^n$ given by:
$$
(z_1,\ldots  ,z_n)\mapsto ((et)^{a_{11}}z_1, (et)^{a_{12}}z_2, \ldots ,(et)^{a_{1n}}z_n)
$$
 which conserves the arguments. Assume now that $A\subset co\mathscr{A}_{V_f}\cap \Critv (\Arg )$,
so  when $t$ is so close to zero then 
 the set $\Log (\Arg^{-1}(A)\cap V_f)$ take place on the two sides of the hyperplane $E^*_{\alpha\gamma}$ in $\Gamma$ dual to $E_{\alpha\gamma}$, because it is the case for the truncation $V_{f^{\Delta_1}}$ which approximate our hypersurface when $t$ tends to zero. So, if one  chooses a coefficients $d_{\alpha}$ and $d_{\gamma}$ such that
 the holomorphic annulus  $\mathscr{Y}$   of equation $d_{\alpha}z^{\alpha} +d_{\gamma}z^{\gamma} = 0$ has   the hyperplane containing $E^*_{\alpha\gamma}$ as its amoeba, and  the hyperplane $H_{\alpha\gamma}$ as its coamoeba, then $V_f\cap \mathscr{Y}$  is nonempty. Let $z_0$ be a point in $ V_f\cap \mathscr{Y}$, hence $\phi_t^{-1}(z_0)\in \phi_t^{-1}(V_f)\cap \phi_t^{-1}(\mathscr{Y})$ and then $\Arg (z_0)\in co\mathscr{A}_{V_{\infty ,\, f^{\Delta_1}}}$. It contradicts Lemma 5.5 if the hyperplane $H_{\alpha\gamma}$ is inner, and  Proposition 5.6 if $H_{\alpha\gamma}$ is external, and then $A$ is contained in the interior of a regular part of the coamoeba or it is discrete.
\end{prooof}

\vspace{0.1cm}

Let $t$ be  a strictly positive real number in $\in ]0;\frac{1}{e}]$, and $H_t$ be the following self
 diffeomorphism of $(\mathbb{C}^*)^n$:
\[
\begin{array}{ccccl}
H_t&:&(\mathbb{C}^*)^n&\longrightarrow&(\mathbb{C}^*)^n\\
&&(z_1,\ldots ,z_n)&\longmapsto&(\parallel z_1\parallel^{-\frac{1}{\log t}}\frac{z_1}{\parallel
  z_1\parallel},\ldots ,\parallel z_n\parallel^{-\frac{1}{\log t}}\frac{z_n}{\parallel z_n\parallel} ).
\end{array}
\]
which defines a new complex structure on $(\mathbb{C}^*)^n$
denoted by $J_t = (dH_t)\circ J\circ (dH_t)^{-1}$ where $J$ is the
standard complex structure.
A $J_t$-holomorphic hypersurface $V_t$ is a hypersurface
holomorphic with respect to the $J_t$ complex structure on
$(\mathbb{C}^*)^n$. It is equivalent to say that $V_t = H_t(V)$ where
$V\subset (\mathbb{C}^*)^n$ is an holomorphic hypersurface with respect to the
standard complex structure $J$ on $(\mathbb{C}^*)^n$.

Recall that the Hausdorff distance between two closed subsets $A,
B$ of a metric space $(E, d)$ is defined by:
$$
d_{\mathcal{H}}(A,B) = \max \{  \sup_{a\in A}d(a,B),\sup_{b\in B}d(A,b)\}.
$$
Here we take $E =\mathbb{R}^n\times (S^1)^n$, with the distance
defined as the product of the
Euclidean metric on $\mathbb{R}^n$ and the flat metric on $(S^1)^n$.

\noindent A complex tropical hypersurface can be defined as follows (see \cite{M1-02} and \cite{M2-04}).
\begin{Def} A complex tropical hypersurface $V_{\infty}\subset
  (\mathbb{C}^*)^n$ is the limit  when
  $t$ tends to zero of a  sequence of a
  $J_t$-holomorphic hypersurfaces $V_t\subset (\mathbb{C}^*)^n$
(with respect to the Hausdorff
  metric on compact sets in $(\mathbb{C}^*)^n$).
\end{Def}

\vspace{0.1cm}

\subsection{Coamoebas of maximally sparse hypersurfaces} ${}$

\vspace{0.2cm}

\noindent Let $V_f\subset (\mathbb{C}^*)^n$ be a hypersurface \, defined by a maximally sparse polynomial 
$\displaystyle{f(z)=}$\\
$\displaystyle{ \sum_{\alpha\in \Verte (\Delta_f)}a_{\alpha}z^{\alpha}}$ (recall that a polynomial $f$ is maximally sparse means that $\supp (f) = \Verte (\Delta_f)$). Let $f_t$ be the family of polynomials defined by \\
$\displaystyle{f_t(z) =\sum_{\alpha\in\Verte (\Delta_f)} a_{\alpha}(et)^{-\Log (a_{\alpha})}z^{\alpha}}$ and $V_t$ their zero locus. We denote by $\displaystyle{V_{\infty , \, f} = \lim_{t\rightarrow 0}H_t(V_t)}$ with respect to the Hausdorff metric on compact sets  of $(\mathbb{C}^*)^n$.

\begin{The} 
With the above notations and assumptions, 
the deformation of $V_f$ given by the family of polynomials $f_t$ satisfies the following:
the coamoeba of the complex tropical hypersurface $V_{\infty ,\, f}$ has the same topology of 
the coamoeba $co\mathscr{A}_{f}$ of $V_f$ (i.e., they are homeomorphic).
\end{The}

\begin{prooof} We will prove that the deformation given by $\{ f_t\}$ defines a bijection between the complement components of the $V_f$'s coamoeba and the complement components of the $V_{\infty ,\, f}$'s coamoeba.  More precisely, we prove that such  deformation conserve the complement components of the coamoeba and thus its topology. Assume that a complement component of the coamoeba is created (resp. disappear) for some $t$. Then there is a created (resp. disappear)
 component of the argument map critical values boundary,
 it means that some edge of the subdivision $\tau_f$ dual to the spine of the amoeba $\mathscr{A}_f$ disappears (resp. created), but it cannot occur because the polynomials are maximally sparse, and thus, the spines of the amoebas $\mathscr{A}_{V_{f_t}}$ are of the same combinatorial type. It remains to show that two different complement components of the $V_f$'s coamoeba
cannot be deformed to the same complement component of  the $V_{\infty ,\, f}$'s coamoeba. Assume on the contrary that there is two complement components $\mathcal{C}_1$ and $\mathcal{C}_2$ of the $V_f$'s coamoeba which are deformed to one complement component $\mathcal{C}$ of the $V_{\infty ,\, f}$'s coamoeba. It  means that one of these two components disappears or the component $\mathcal{C}$ is not convex, and then  we have a contradiction in both  cases.

\end{prooof}

\vspace{0.2cm}

\section{Examples of complex algebraic plane curves coamoebas}

\begin{itemize}
\item[(1)]\,
Let $V_{f_{\lambda}}$ be the curve in $(\mathbb{C}^*)^2$ defined by the following polynomial:
$$
f_{\lambda}(z,w) = w^2-\lambda w + 2zw -z^2w +1.
$$
Let $f_{\lambda ,\, 1}(z,w)=  w^2-\lambda w + 2zw -z^2w $, so $V_{f_{\lambda ,\, 1}}$ is just the parabola of example 1.
Let $f_{\lambda ,\, 2}(z,w)=  -\lambda w + 2zw -z^2w +1$, hence $V_{f_{\lambda ,\, 2}}$ is the set of points
$(z,w)\in (\mathbb{C}^*)^2$ such that :
$$
w=\frac{1}{z^2-2z+\lambda}.
$$
This means that $\arg (w_2) = - \arg (w_1)$ mod $2\pi$. Hence the coamoeba of the curve defined by $f_{\lambda}$ is as in the figure 8 on the left.
\begin{figure}[h!]
\begin{center}
\includegraphics[angle=0,width=0.3\textwidth]{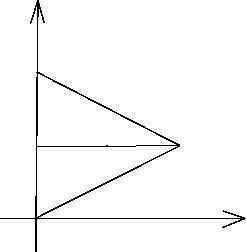}
\caption{The Newton Polygon of example (1) and its subdivision.}
\label{c}
\end{center}
\end{figure}

\begin{figure}[h!]
\begin{center}
\includegraphics[angle=0,width=0.3\textwidth]{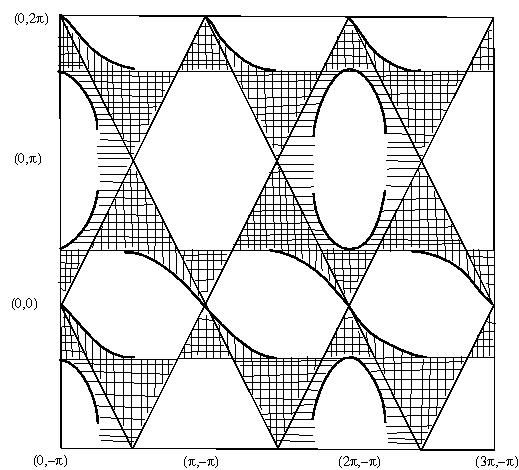}
\includegraphics[angle=0,width=0.3\textwidth]{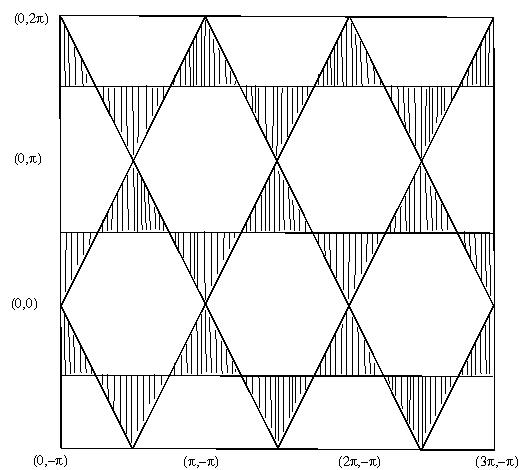}
\caption{Example (1): on the left the coamoeba when $\lambda\ne 0$ and the curve is non-Harnack, and the coamoeba when $\lambda = 0$ and the triangulation is trivial  on the right.}
\label{c}
\end{center}
\end{figure}


\item[(2)]\,
Let $V_{f_{\lambda}}$ be the curve in $(\mathbb{C}^*)^2$ defined by the following polynomial:
$$
f_{\lambda}(z,w) = zw^2+z^2w+z+w+ \lambda zw.
$$
Let $f_{\lambda ,\, 1}(z,w)= zw^2+z+w+ \lambda zw$. Hence $V_{f_{\lambda ,\, 1}}$ is just
a reparametrization of the parabola of example 1. We can see that $z = -\frac{w}{1+w^2+\lambda w}$.\\
Let $f_{\lambda ,\, 2}(z,w)= zw^2+z^2w+z+ \lambda zw = z(1+zw+w^2+\lambda w) $, hence
$V_{f_{\lambda ,\, 2}}$ is the set of points $(z,w)\in (\mathbb{C}^*)^2$ such that :
$$
z=-\frac{1+w^2+\lambda w}{w}.
$$
It means that $\arg (z_2) = - \arg (z_1)$ mod $2\pi$, where $z_1$ (resp. $z_2$)
denotes  the first coordinate of a point in $V_{f_{\lambda ,\, 1}}$ (resp. in
$V_{f_{\lambda ,\, 2}}$) . As in example 1, we have the figures 9 on the top right.

\begin{figure}[h!]
\begin{center}
\includegraphics[angle=0,width=0.3\textwidth]{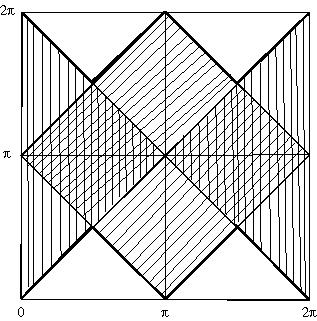}
\includegraphics[angle=0,width=0.5\textwidth]{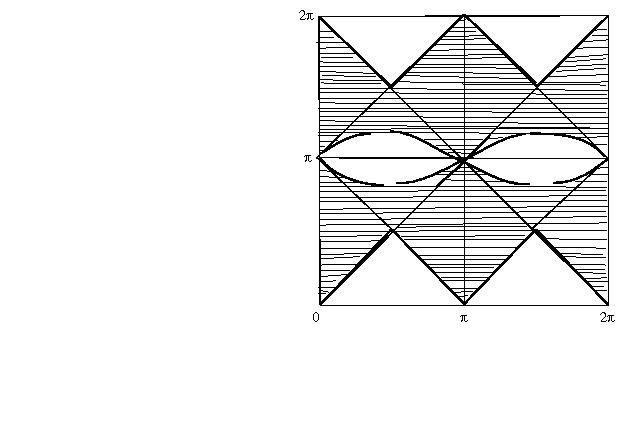}
\includegraphics[angle=0,width=0.3\textwidth]{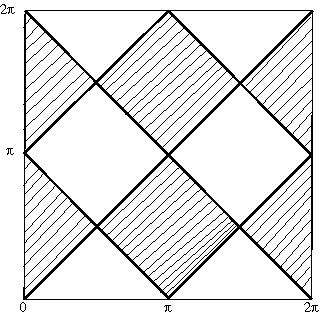}
\caption{Coamoeba of example (2) in three cases, the first coamoeba is the one of  the Harnack case, the second coamoeba of a non-Harnack case and  the coefficient $\lambda\ne 0$, and the last case is the case when $\lambda = 0$ and the subdivision is trivial.}
\label{c}
\end{center}
\end{figure}

\item[(3)]\,
We give an example of a Newton polygon $\Delta$ that defines not any 
real curve 
  with maximal number of coamoeba complement components, but 
 this maximal number is realized by a complex curve. Let $\Delta$ be the polygon with vertices $(1;0),\,
 (0;1),\, (1;2)$, and $(3;1)$ (see figure 10 for the polygon and its subdivision dual to the spine of the
 amoeba). In this case we prove that no real polynomial can realize the maximal number of coamoeba complement components (the maximal number in the real case is five, and the coamoeba is given in figure 11 on the left for some
 real coefficients ), but the complex curve defined by the  complex polynomial $f(z,w) = 
e^{i\alpha}w+z+zw^2+z^3w$ with $0<\alpha <\pi$, has a coamoeba with maximal number of complement components 
(i.e. six components, see figure 11 on the right).

\begin{figure}[h!]
\begin{center}
\includegraphics[angle=0,width=0.3\textwidth]{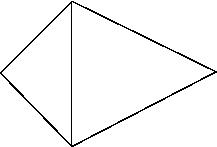}\qquad\qquad\qquad
\includegraphics[angle=0,width=0.3\textwidth]{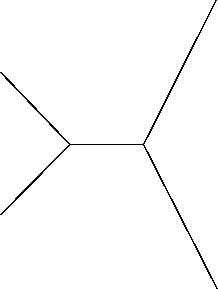}
\caption{The subdivision of the Newton polygon and its dual}
\label{c}
\end{center}
\end{figure}

\begin{figure}[h!]
\begin{center}
\includegraphics[angle=0,width=0.3\textwidth]{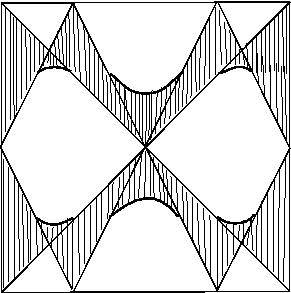}\qquad\qquad\qquad
\includegraphics[angle=0,width=0.4\textwidth]{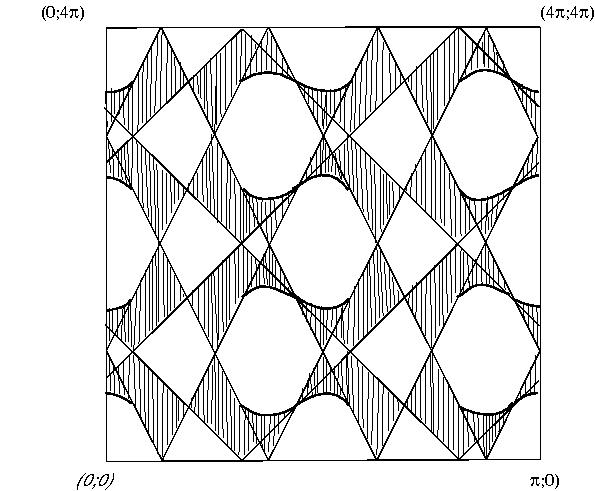}
\caption{Example (3): on the left the coamoeba of a real curve with the same  maximal number of complement components (i.e., 5 components; the picture here is in one fundamental domain) and on the right the
coamoeba of a complex curve with a maximal number of complement components (6 components; the picture here is in four fundamental domains).}
\label{c}
\end{center}
\end{figure}

\end{itemize}

\end{document}